\documentclass[12pt]{amsart}
\usepackage{amsmath}
\usepackage{amsthm}
\usepackage{latexsym}
\usepackage{amssymb}

\newtheorem{theorem}{Theorem}
\newtheorem{proposition}{Proposition}

\newtheorem{definition}{Definition}

\newtheorem{remark}{Remark}
\newcounter{obsctr}

\renewcommand{\theequation}{\thesection.\arabic{equation}}

\begin{document}
\def\A {{\mathcal{A}}}
\def\D {{\mathcal{D}}}
\def\R {{\mathbb{R}}}
\def\N {{\mathbb{N}}}
\def\C {{\mathbb{C}}}
\def\Z {{\mathbb{Z}}}
\def\l {\ell}
\def\tm {\tilde{m}}
\def\ml {multline}
\def\multiline {\multline}
\def\lessim {\lesssim}
\def\phi{\varphi}
\def\epsilon{\varepsilon}
\def\olm{\overline{L_m}}
\def\ol{\overline{L}}
\def\oz{\overline{z}}
\def\be{\begin{equation}}
\def\\[{\begin{equation}}
\def\ee{\end{equation}}
\def\\]{\end{equation}}
\title{Analyticity for singular sums of squares of
degenerate vector fields}    
\author{David S. Tartakoff}
\address{Department of Mathematics, University
of Illinois at Chicago, m/c 249, 851 S.
Morgan St., Chicago IL  60607, USA}
\email{dst@uic.edu}
\date{\today}
\begin{abstract} Recently, J.J. Kohn
in \cite{K2005} proved hypoellipticity for
\begin{equation} P=L\overline{L} +
\overline{L}|z|^{2k}L \hbox{ \quad with \quad}
L={\partial
\over
\partial z} + i\overline{z}{\partial \over
\partial t},\tag{$*_k$}\end{equation}
i.e.,
$$-P=\overline{L}^*\overline{L} +
(\overline{z}^kL)^*\overline{z}^kL,$$ 
 a singular sum
of squares of complex vector fields on the
complex Heisenberg group, an operator which exhibits a
{\it loss} of
${k-1}$ derivatives. Subsequently, in
\cite{DT2005}, M. Derridj and D. S. Tartakoff
proved analytic hypoellipticity for this operator using
rather different methods going back to \cite{Ta1978},
\cite{Ta1980}. Together with A. Bove and J. J. Kohn,
Derridj and Tartakoff in \cite{BDKT2005} gave an
alternative proof of the
$C^\infty$ hypoellipticity of $*_k$ in the style of
\cite{DT2005}. In this note, we consider the
equation 
\begin{equation} P=L_m\overline{L_m} +
\overline{L_m}\,|z|^{2k}L_m \hbox{ \; with \;}
L_m={\partial
\over
\partial z} + i\overline{z}|z|^{2m}{\partial \over
\partial t},\tag{$*_{m,k}$}\end{equation}
for which the underlying manifold is only of finite
type, and prove analytic hypoellipticity. This operator
is subelliptic with large loss of derivatives but the
exact loss is irrelevant for analytic hypoellipticity.
\end{abstract}
\maketitle
\pagestyle{myheadings}
\markboth{D.S. Tartakoff}
{Analyticity for degenerate vector fields}
\section{Introduction and statement of theorems}
\renewcommand{\theequation}{\thesection.\arabic{equation}}
\setcounter{equation}{0}
\setcounter{theorem}{0}
\setcounter{proposition}{0}  
\setcounter{lemma}{0}
\setcounter{corollary}{0} 
\setcounter{definition}{0}
\setcounter{remark}{0}

J.J. Kohn's recent paper \cite{K2005}, inspired by work
of Siu on singular metrics and the implied applications,
studied a singular sum of squares of complex vector
fields on the (complex) Heisenberg group: 
\begin{equation} P=-L\overline{L} -
\overline{L}|z|^{2k}L \hbox{ \quad with \quad}
L={\partial
\over
\partial z} + i\overline{z}{\partial \over
\partial t},\tag{$*_k$}\end{equation}
i.e.,
$$P=\overline{L}^*\overline{L} +
(\overline{z}^kL)^*\overline{z}^kL,$$ 
He showed that this operator was hypoelliptic but loses
$k-1$ derivatives (in Sobolev norms), and a note by
Christ showed that the addition of the square of
another vector field, ${\partial^2\over \partial
s^2}$ to $P$, destroyed hypoellipticity
completely.  

Subsequently, in
\cite{DT2005}, M. Derridj and D. S. Tartakoff
proved analytic hypoellipticity for this operator using
rather different methods going back to \cite{Ta1978},
\cite{Ta1980}. In 
\cite{BDKT2005}, together with
A. Bove and J. J. Kohn, M. Derridj and D.S. Tartakoff 
then gave an alternative proof of the
 $C^\infty$ hypoellipticity of $*_k$ in the style of
\cite{DT2005}. 

In this paper, we consider the more degenerate 
equation 
\begin{equation} P= P_m=L_m\overline{L_m} +
\overline{L_m}\,|z|^{2k}L_m \hbox{ \; with \;}
L_m={\partial
\over
\partial z} + i\overline{z}|z|^{2m}{\partial \over
\partial t},\tag{$*_{m,k}$}\end{equation}
for which the underlying manifold is only of finite
type, and prove analytic
hypoellipticity for this operator. 
\begin{theorem} The operator $P_m$ is analytic
hypoelliptic: $P_mu\in C^\omega \rightarrow u\in
C^\omega.$

\end{theorem}

\section{The {\it a priori} estimate}
\renewcommand{\theequation}{\thesection.\arabic{equation}}
\setcounter{equation}{0}
\setcounter{theorem}{0}
\setcounter{proposition}{0}  
\setcounter{lemma}{0}
\setcounter{corollary}{0} 
\setcounter{definition}{0}
\setcounter{remark}{0}

As in \cite{DT2005} we work with $v\in C_0^\infty$ and
denote by
$\Lambda$ the `pseudodifferential' operator with symbol
$\lambda (z,t; \zeta,\tau) = (1+|\tau|^2)^{1/2}.$ While
$\Lambda$ is not a true pseudodifferential operator,
for $\zeta$ different from zero the operator $P_m$ is
elliptic, hence analytic hypoelliptic, near $z=0$ and
for $z$ different from $0,$ the operator
$P_m$ has symplectic characteristic variety and satisfies
maximal estimates in $L$ and $\overline{L},$ hence is
even analytic hypoelliptic there by \cite{Ta1978},
\cite{Tr1978}, \cite{Ta1980}. Thus we will
need only study the behavior of powers of
$\partial / \partial t$ applied to the solution $u$
locally. We shall do so in $L^2$ norms with appropriate
powers of $\Lambda.$

The {\it a priori} estimate satisfied by $P_m$ is
actually better than that satisfied by $P_0$ but as is
well known, analyticity does not use the precise degree
of subellipticity in an important way. Thus we give a
simple proof of the estimate given in Kohn's paper
\cite{K2005}.

\begin{proposition} For all values of $m,$ we have
$$\|\olm v\|^2 +
\|\oz^kLv\|^2 + \|v\|^2_{-\frac{k-1}{2}}\leq C
|(P_mv,v)|,
\quad v\in C_0^\infty.$$
\end{proposition}

For any $r<0,$ since 
$\Lambda$ commutes with everything, we have:
$$\|\Lambda^r v\|^2 = ([L_m,z]\Lambda^rv,
\Lambda^rv)=  (z\Lambda^rv,\Lambda^r\overline{L_m}v) 
- (L_m\Lambda^{r-\frac{1}{2}}v,
\overline{z}\Lambda^{r+\frac{1}{2}}v)$$ 
so 
\begin{equation}\label{vr}\|\Lambda^rv\|^2 \leq
s.c.\|\Lambda^rv\|^2 + l.c.\,\|\olm v\|^2 + 
s.c.\|L_m\Lambda^{r-\frac{1}{2}}v\|^2 +
l.c.\|\overline{z}\Lambda^{r+\frac{1}{2}}v\|^2\end{equation}

Now integrating by parts and commuting $L_m$ and
$\olm$ we have 
$$\|L_m\Lambda^{r-\frac{1}{2}}v\|^2 \leq 
\|\olm\Lambda^{r-\frac{1}{2}}v\|^2 +
\|\Lambda^{r}v\|^2$$
so that equation (\ref{vr}) becomes
\begin{equation}\label{vr2}\|\Lambda^rv\|^2
\lesssim\|\olm v\|^2 + 
\|\overline{z}\Lambda^{r+\frac{1}{2}}v\|^2
\lesssim |(P_m v, v)| +
\|\overline{z}\Lambda^{r+\frac{1}{2}}v\|^2\end{equation}

In the last term we integrate by parts, and obtain
$$\|\overline{z}\Lambda^{r+\frac{1}{2}}v\|^2 = 
(\overline{z}\Lambda^{r+\frac{1}{2}}v,
\overline{z}\Lambda^{r+\frac{1}{2}}v) =
(z\oz\Lambda^{r+1}v, \Lambda^r v)|$$
$$ \leq
s.c.\|\Lambda^rv\|^2 + \l.c.\|z^2\Lambda^{r+1}v\|^2$$
so that, iterating and absorbing the first term on the
left in (\ref{vr2}),
\begin{equation}\label{vr3}
\|\Lambda^r v\|^2 \lesssim |(P_mv,v)| +
\|z^{k-1}\Lambda^{r+\frac{k-1}{2}}v\|^2
\end{equation}

To finish the
derivation of the {\it a priori} estimate, we write
$kz^{k-1}=[L_m,z^k],$ so that
\begin{equation}\label{est:zk-1} k^2\|z^{k-1}w\|^2 =
([L_m,z^k] w, z^{k-1}w) \end{equation}
$$ = (L_mz^k w,
z^{k-1}w) - (z^kL_m w,
z^{k-1}w)= \tilde{A} +
\tilde{B}.$$
Then  
$$|\tilde{A}| \lesssim |(z^{k-1} w,
\overline{z}\overline{L_m}z^{k-1}w)| \lesssim s.c.
\|z^{k-1}w\|^2 +
l.c.\|\overline{L_m}w\|^2$$
since $[\overline{L_m}, z^{k-1}]=0,$ while 
$$|\tilde{B}|\sim s.c. \|z^{k-1}w\|^2 + l.c.
\|z^kL_mw\|^2.$$ The terms with the small constant
($s.c.$) will be absorbed on
the left hand side of (\ref{est:zk-1}), yielding:
\begin{equation}\label{zk-1vr} \|z^{k-1}w\|^2 \leq
C\{\| \overline{L_m}w\|^2 + \|
\oz^kL_mw\|^2\}\end{equation}

Combining this estimate with (\ref{vr}), since
$\|z^kw\|=\|\overline{z}^kw\|,$ leads to
\begin{equation}\label{est:ap1}\|\Lambda^{{-{k-1\over
2}}}v\|^2 \leq C\{\|\overline{L_m}v\|^2
 + \|\overline z^kL_mv\|^2\} =
C|(P_mv,v)|. \end{equation}  
with
no errors (except the `microlocalizing ones, i.e.,
staying out of the elliptic region and the region $z$
not equal to $0$, where $P_m$ has symplectic
characteristic variety and
\cite{Ta1978},
\cite{Tr1978}, and \cite{Ta1980} apply). 
\begin{remark} If one wishes to use the full
pseudodifferential operator with symbol
$(1+|\zeta|^2+|\tau|^2)^{1/2}$ instead of
$\Lambda^{1/2},$ the derivation is virtually identical,
with the addition of a lower order norm on the
right of (\ref{est:ap1}).
\end{remark}

\section{The localization of powers of $T$ and its
commutators}
\renewcommand{\theequation}{\thesection.\arabic{equation}}
\setcounter{equation}{0}
\setcounter{theorem}{0}
\setcounter{proposition}{0}  
\setcounter{lemma}{0}
\setcounter{corollary}{0} 
\setcounter{definition}{0}
\setcounter{remark}{0}

To apply (\ref{est:ap1}) we must replace $v$ by a
localization of powers of ${\partial/\partial t},$ or,
for convenience, we will localize powers of 
$$T=-\frac{2i}{m+1}{\partial\over\partial t}$$  applied
to the solution
$u$. Note that 
$$[L_m, \overline{L_m}\,] = (m+1)^2|z|^{2m}T,$$
which will be of technical use later. 

For now we treat the solution
$u$ as being smooth, and concentrate on the estimates,
remarking that by multiplying the Fourier Transform 
$\tilde{u}(z,\oz,\tau)$ of $u$ in $t$ by dilations of
$\Psi(\tau)$ equal to zero for
$|\tau|>2$ but
identically equal to one for 
$|\tau|<1,$ once we obtain uniform estimates (in
$\epsilon$) for derivatives of the
inverse transform of $\Psi(\tau/\epsilon)(\tilde{u})$ we
may let $\epsilon\rightarrow 0$ to obtain the same
estimates for the same derivatives of $u.$ 

The localization in space may be taken independent of
$z,\oz$ since taking a product of two localizers,
$\phi_1(t)\phi_2(z,\oz)$ in terms where no derivatives
land on $\phi_2$ it just survives to the left of
$P_m=f$ while in terms with derivatives, those terms
are supported in the region where $z$ is different from
zero, i.e., in regions where the result is known, as
remarked above. 

But even with $\phi$ depending on $t$ alone,
the localization must be done very carefully. For
example, the first 
 bracket with $P_m$ which we encounter will contain 
$(L_m\phi (t)) T^p \sim \overline{z}|z|^{2m}\phi'T^p,$
which is problematic for any value of $m.$ 

In Kohn's
work
\cite{K2005} ($m=0$), the 
$\oz$ (or $z$) in front of each derivative of $\phi$ is
carefully followed, and shown to provide, after
some work, a gain of $1/2$ derivative. Analyticity
was not considered in that paper, nor does it seem
likely that it could be shown by those methods. 

Derridj and Tartakoff found in
\cite{DT2005} that an entirely different
approach, involving a delicately balanced localization
of
$T^p,$ led to analyticity rather directly. 

Here, to handle the case $m>0,$ which corresponds to a
CR manifold of finite type, we find it simpler
to use a somewhat differently balanced localization.
And with this localization both $C^\omega$ and
$C^\infty$ hypoellipticity come together. 

\subsection{The case $m=0$}

For motivation, however, we first we give the definition
in the case $m=0.$
\begin{definition} 
For any pair of
non-negative integers, $(p_1,p_2),$ let
$$(T_0^{p_1,p_2})_\phi = \sum_{{a\leq p_1} 
\atop {b\leq p_2}}{\overline{L_m}^a\circ 
\overline{z}\,^a\circ 
T\,^{p_1-a}\circ
\phi^{(a+b)}\circ T\,^{p_2-b}\circ
{z}^b\circ L_m^b
\over a!b!},$$
where $$\phi^{(r)} =\left(-i\frac{\partial}{\partial
t}\right)^r\phi.$$
\end{definition}

\noindent Note that
the leading term (with
$a+b=0)$ is merely
$T^{p_1}\circ \phi \circ T^{p_2}$ which is equal to the
operator 
$T^{p_1+p_2}$ on any open set $\Omega_0$ where $\phi
\equiv 1.$

We have the extraordinary commutation relations: 

\begin{proposition}\label{3.1} 
$$[L_0, (T_0^{p_1,p_2})_\phi] \equiv
(T_0^{p_1,p_2-1})_{\phi'}\circ L_0,$$
$$ [\overline{L_0},
(T_0^{p_1,p_2})_\phi] \equiv 
\overline{L_0}\circ (T_0^{p_1-1,p_2})_{\phi'},
$$
$$[(T_0^{p_1,p_2})_\phi,z] =
z\circ (T_0^{p_1,p_2-1})_{\phi'},
$$
and
$$
[(T_0^{p_1,p_2})_\phi,\overline{z}] =
(T_0^{p_1-1,p_2})_{\phi'}\circ \overline{z},
$$
\vskip.1in\noindent
where the $\equiv$ denotes modulo
$C^{p_1-p_1'+p_2-p_2'}$ terms of the form
$${\overline{L_0}^{p_1-p_1'}\circ
\overline{z}^{p_1-p_1'}\circ 
T^{p_1'}\circ\phi^{(p_1'+p_2'+1)}\circ T^{p_2'}
\circ{z}^{p_2-p_2'}\circ
{L_0}^{p_2-p_2'}\over (p_1-p_1')!(p_2-p_2')!}$$
with either $p_1'=0$ or
$p_2'=0,$ i.e., terms where all free $T$
derivatives have been eliminated on one side of
$\phi$ or the other. 
\end{proposition}
\begin{proof}
The proof is a straightforward calculation involving a
shift of index in the definition of
$(T_0^{p_1,p_2})_\phi.$
\end{proof}

\subsection{The case $m>0$}

For $m>0,$ we use a related, but somewhat
different definition of the
localization, owing to multiple brackets. Parts of this
discussion are in the style of \cite{DT1988},
\cite{BT2005}. 
\begin{definition} For $m>0,$ and $(p_1, p_2)$ as
 above, set 
$$\phi^{(d)} = 
\left(i\frac{\partial}{\partial t}\right)^{d}\phi(t),$$
$$N_b = \sum_{b'leq b}
A^b_{b'}{(zL_m)^{b'}\over b'!},$$ where the $A^b_{b'}$
(real) are to be determined subject to $A^b_b=1,$ and 
 $$\tilde{N}_a = \sum_{a'\leq a}
A^a_{a'}{(-\olm\oz)^{a'}\over
a'!}=(N_a)^*$$
with the same coefficients, and finally we set
$$(T_m^{p_1,p_2})_\phi = \sum_{{a\leq p_1} 
\atop {b\leq p_2}}\tilde{N}_a\circ  
T\,^{p_1-a}\circ
\phi^{(a+b)}\circ T\,^{p_2-b}\circ
N_b.$$

\end{definition}

We have
$$[\olm, (T_m^{p_1,p_2})_\phi] = [\olm, \sum_{a\leq
p_1,
b\leq p_2}\tilde{N}_a  T^{p_1-a}
\phi^{(a+b)} T^{p_2-b} N_b]$$
\begin{equation}\label{LTpphi}
=\sum_{a\leq
p_1,b\leq p_2}\bigg\{[\olm,
\tilde{N}_a] T^{p_1-a}
\phi^{(a+b)} T^{p_2-b} N_b
\end{equation}
$$-\tilde{N}_a T^{p_1-a} z|z|^{2m}
\phi^{(a+b+1)} T^{p_2-b} N_b
+ \tilde{N}_a  T^{p_1-a}
\phi^{(a+b)}
T^{p_2-b}[\olm, N_b]\bigg\}.$$

The last two terms on the right must cancel, to
preserve the balance, since both disturb the balance
between derivatives on $\phi$ and gain in powers of $T.$
We will choose the coefficients $A^b_{b'}$ of
${N}_b$ in such a way that, modulo
acceptable errors,   
\begin{equation}\label{recursion} [\olm,N_b] =
z|z|^{2m}T{N}_{b -1}.
\end{equation}
This will provide the needed cancellation via a shift
of index in $a$ in the sum just as in the case with
$m=0.$ The corresponding relation for brackets with
$L_m$ will follow by taking adjoints: again modulo
acceptable errors, 
\begin{equation}\label{recursion2}[L_m,\tilde{N}_a] = 
-\tilde{N}_{a -1}\oz|z|^{2m}T.
\end{equation}

Condition (\ref{recursion}) reads, using the
definition of ${N}_b,$ reads:
$$\sum_{b'=0}^b
A^b_{b'}{1\over b'!}\left[\olm, (zL_m)^{b'}\right] =
-\oz|z|^{2m}T\sum_{b'=0}^{b -1} A^{b
-1}_{b'}{(zL_m)^{b'}\over b'!}.$$

Expanding the brackets and keeping all factors of
$zL_m$ to the right, 
$$ {1\over b'!}\left[\olm, (zL_m)^{b'}\right]={1\over
b'!}\sum_{1\leq b''\leq b'}{b'\choose
b''}ad_{zL_m}^{b''}(\olm)(zL_m)^{b'-b''}$$
$$={1\over
b'!}\sum_{1\leq b''\leq b'}{b'\choose
b''}ad_{zL_m}^{b''-1}((m+1)^2z|z|^{2m}T)(zL_m)^{b'-b''}
$$
$$= (z|z|^{2m}T)\sum_{1\leq b''\leq
b'}\frac{(m+1)^{b''+1}}{b''!}
\,\frac{(zL_m)^{b'-b''}}{(b'-b'')!}.$$
The condition (\ref{recursion}) thus requires, renaming
$b'-b''$ as $\overline{b}$ on the right just above,
\begin{equation}\label{rrB}\sum_{b''=1}
^{b-\tilde b} A^b_{\tilde b
+ b''}\frac{(m+1)^{b''+1}}{b''!}
=A^{b-1}_{\tilde{b}}.
\end{equation}

Fortunately, we have investigated these
equations in \cite{DT1988} and, citing a result in the
book by Hirzebruch \cite{H1966} have explicit
solutions $\tilde{A}^*_*$ in the case that $m=0$, unique
under the conditions that $\tilde{A}^q_0 = (-1)^q,$
namely 
$$\tilde{A}^r_s=\left(\left(\frac{t}
{e^t-1}\right)^{r+1}\right)^{(r-s)}(0)/(r-s)! $$
and for positive $m$ we merely dilate by the factor of
$m+1:$ 
$$A^r_s=\left(\left(\frac{\frac{t}{m+1}}
{e^{\frac{t}{m+1}}-1}
\right)^{r+1}\right)^{(r-s)}(0)/(r-s)!$$ 

In addition, we will also need good expressions for
the other brackets: we compute  
$$[L_m, N_b] = [L_m,\sum_{b'=0}^b
A^b_{b'}\frac{(zL_m)^{b'}}{b'!}] 
= \sum_{{1\leq b''\leq b'}\atop {b'\leq b}}
A^b_{b'}\,\frac{1}{b''!}\,
\frac{(zL)^{b'-b''}}{(b'-b'')!}\circ
L_m$$
$$
[N_b, z] = z\circ \sum_{{1\leq b''\leq b'}\atop {b'\leq
b}} A^b_{b'}\,\frac{1}{b''!}\,
\frac{(zL)^{b'-b''}}{(b'-b'')!}
$$
\begin{equation}\label{LoLN}
[\olm,
\tilde{N}_a] =
-\,\olm \circ \sum_{{1\leq a''\leq a'}\atop {a'\leq a}}
A^a_{a'}\,\frac{1}{a''!}\,
\frac{(-\olm\oz)^{a'-a''}}{(a'-a'')!},
\end{equation}
$$[\oz,\tilde{N}_a] =
\sum_{{1\leq a''\leq a'}\atop {a'\leq a}}
A^a_{a'}\,\frac{1}{a''!}\,
\frac{(-\olm\oz)^{a'-a''}}{(a'-a'')!}\circ \oz.
$$

In order to recognize these sums as $N$'s or
$\tilde{N}$'s, we need to be able to shift the lower
indices on $A^a_{a'}$
down by one. But this also we have done in
\cite{DT1988}, with the result that 
\begin{proposition} For any $r,s,$ and $c,$ 
$$A^r_s=\sum_{j=0}^{r-s}S_j^{r-s}\,A^{r-(c+j)}_{s-c}
$$
where
$$|S^k_\ell|\leq C^k.$$
\end{proposition}
These brackets, then, together with the Proposition,
immediately translate, setting $b''=c$ and
$\tilde{b}=b'-b'',$ into:
$$[L_m, N_b] = \sum_{\tilde{b}\leq b-c}
\frac{1}{c!}A^b_{\tilde{b}+c}\,\,
\frac{(zL)^{\tilde{b}}}{\tilde{b}!}\circ
L_m=
\,
\sum_{\tilde{b}\leq
b-c-j}\frac{1}{c!}S^{c+j}_j\,A^{b-c-j}_{\tilde{b}}\frac{(zL)^{\tilde{b}}}{\tilde{b}!}\circ
L_m$$
or 
$$[L_m, N_b] = \sum_{{c+j\leq
b}\atop 1\leq c}\frac{1}{c!}S^{c+j}_j N_{b-c-j}\circ
L_m$$

Similarly, 
$$[\olm, \tilde{N}_a] = -\olm\circ \sum_{{c+j\leq
a}\atop 1\leq c}\frac{1}{c!}S^{c+j}_j
\tilde{N}_{a-c-j},$$ and 
$$[ \oz, \tilde{N}_a] = \sum_{{c+j\leq
a}\atop 1\leq c}\frac{1}{c!}S^{c+j}_j
\tilde{N}_{a-c-j}\circ \oz,$$ and
$$[{N}_b, z] = z\circ \sum_{{c+j\leq
b}\atop 1\leq c}\frac{1}{c!}S^{c+j}_j \tilde{N}_{b-c-j},$$

These wonderful commutation relations mean that the
whole localization $(T^{p_1,p_2})_\phi$ may be commuted
meaningfully with the vector fields $L_m, \olm$ and
with $z, \oz:$

\begin{proposition}\label{3.3} Modulo terms in which
either
$p_1$ or $p_2$ has been reduced to zero, and in view of
the cancellations ensured by (\ref{recursion}), 
\begin{equation}\label{LmTpphi}[L_m,
(T_m^{p_1,p_2})_\phi]\equiv
\sum_{{1\leq c, 0\leq j}\atop c+ j\leq
p_2}\frac{1}{c!}S^{c+j}_j
(T_m^{p_1,p_2-(c+j)})_{\phi^{(c+j)}}\circ
L_m\end{equation}
\end{proposition}
\begin{proof}
$$[L_m, (T_m^{p_1,p_2})_\phi]\equiv \sum_{{a\leq p_1} 
\atop {b\leq p_2}}\tilde{N}_a\circ  
T\,^{p_1-a}\circ
\phi^{(a+b)}\circ T\,^{p_2-b}[L_m, N_b]$$
$$\equiv \sum_{{a\leq p_1} 
\atop {b\leq p_2}}\sum_{{c+j\leq
b}\atop 1\leq c}\frac{1}{c!}S^{c+j}_j\tilde{N}_a\circ  
T\,^{p_1-a}\circ
\phi^{(a+b)}\circ T\,^{p_2-b} N_{b-c-j}\circ L_m$$
$$\equiv \sum_{{a\leq p_1} 
\atop {b\leq p_2}}\sum_{{c+j\leq
b}\atop 1\leq c}\frac{1}{c!}S^{c+j}_j\tilde{N}_a  
T\,^{p_1-a}
{\phi^{(c+j)}}^{(a+b-c-j)} T\,^{p_2-(c+j)-(b-c-j)}
N_{b-c-j} \circ L_m
$$
$$\equiv
\sum_{{1\leq c, 0\leq j}\atop c+ j\leq
p_2}\frac{1}{c!}S^{c+j}_j
(T_m^{p_1,p_2-(c+j)})_{\phi^{(c+j)}}\circ L_m$$
\end{proof}

Similarly we state, and omit the proofs, which are
virtually identical to that of the previous
proposition,  

\begin{proposition} 
\begin{equation}\label{LmbarTpphi} [\olm,
(T_m^{p_1,p_2})_\phi]\equiv -\olm\circ \sum_{{1\leq c,
0\leq j}\atop c+ j\leq p_1}\frac{1}{c!}S^{c+j}_j
(T_m^{p_1-(c+j),p_2})_{\phi^{(c+j)}}, 
\end{equation}
\begin{equation}\label{zbarTpphi}[\oz,
(T_m^{p_1,p_2})_\phi]\equiv
\sum_{{1\leq c, 0\leq j}\atop c+ j\leq
p_1}\frac{1}{c!}S^{c+j}_j
(T_m^{p_1-(c+j),p_1})_{\phi^{(c+j)}}\circ
\oz,\end{equation}
 and
\begin{equation}\label{zTpphi}[z,
(T_m^{p_1,p_2})_\phi]\equiv -\sum_{{1\leq c, 0\leq
j}\atop c+ j\leq p_2}\frac{1}{c!}S^{c+j}_j
(T_m^{p_1,p_2-(c+j)})_{\phi^{(c+j)}}\circ
z\end{equation}

\end{proposition}

What these commutation relations mean is that we may
move the vector fields of $P_m$ past
$(T_m^{p_1,p_2})_\phi$ freely, at each stage incurring
errors with the same vector fields and a gain in
derivatives in $(T_m^{p_1,p_2})_\phi.$ Thus we may
iterate the {\it a priori} inequality modulo errors of
nearly arbitrarily low order - all of the $\equiv$
signs above mean that we will ultimately arrive at
errors where either $p_1=0$ or $p_2=0.$

 So we insert first 
$v=(T_m^{p_1,p_2})_\phi u$ into 
(\ref{est:ap1}), then bring
$(T_m^{p_1,p_2})_\phi$ to the left of
$P_m=L\overline{L}+
 \overline{L}z^k\overline{z}^kL,$ now writing $L$
instead of $L_m,$ since the formal expansions of
the brackets are insensitive to
$m.$ 
so that we have:
$$\|\overline{L}(T_m^{p_1,p_2})_\phi
u\|_0^2 +
\|\overline{z}^k{L}(T_m^{p_1,p_2})_\phi
u\|_0^2 +\|\Lambda^{-{k-1\over
2}}(T_m^{{p\over 2}, {p\over 2}})_\phi
u\|_0^2$$
\begin{equation}\label{est:ap2}\lesssim
|(P_m(T_m^{p_1,p_2})_\phi u, (T_m^{p_1,p_2})_\phi u)_{L^2}|
\end{equation}
$$\lesssim
|((T_m^{p_1,p_2})_\phi P_mu, (T_m^{p_1,p_2})_\phi u)_{L^2}| + |([P_m,(T_m^{p_1,p_2})_\phi] u, (T_m^{p_1,p_2})_\phi u)_{L^2}|$$ and
by the above bracket relations, modulo the
same terms as above where all $T$'s from one side of
$\phi$ or the other have been `converted' into $L$'s or
$\overline{L}$'s, we have
$$([P_m,(T_m^{p_1,p_2})_\phi] u,
(T_m^{p_1,p_2})_\phi u) \equiv$$
$$= ([L\overline{L},(T^{p_1,p_2})_\phi]
u, (T^{p_1,p_2})_\phi u) + ([
\overline L z^k\overline{z}^k{L},(T^{p_1,p_2})_\phi] u, (T^{p_1,p_2})_\phi u)$$
$$= ([L,(T^{p_1,p_2})_\phi]\overline{L}
u, (T^{p_1,p_2})_\phi u) +
(L[\overline{L},(T^{p_1,p_2})_\phi] u,
(T^{p_1,p_2})_\phi u)$$
$$+([
\overline L ,(T^{p_1,p_2})_\phi]z^k\overline{z}^k{L} u, (T^{p_1,p_2})_\phi u)+ (\overline L[
 z^k,(T^{p_1,p_2})_\phi] \overline{z}^k{L}u, (T^{p_1,p_2})_\phi u)$$
$$+(\overline L z^k[
\overline{z}^k,(T^{p_1,p_2})_\phi]{L} u, (T^{p_1,p_2})_\phi u)+ 
(\overline L z^k\overline{z}^k[{L},(T^{p_1,p_2})_\phi] u, (T^{p_1,p_2})_\phi
u)$$
$$\equiv \sum_{{1\leq c, 0\leq j}\atop c+ j\leq
p_2}\frac{1}{c!}S^{c+j}_j(
(T_m^{p_1,p_2-(c+j)})_{\phi^{(c+j)}}
L\overline{L}u, (T_m^{p_1,p_2})_\phi u)$$
$$-\sum_{{1\leq c,
0\leq j}\atop c+ j\leq p_1}\frac{1}{c!}S^{c+j}_j(L\overline{L} 
(T_m^{p_1-(c+j),p_2})_{\phi^{(c+j)}} u, (T_m^{p_1,p_2})_\phi u)$$
$$- \sum_{{1\leq c,
0\leq j}\atop c+ j\leq p_1}\frac{1}{c!}S^{c+j}_j
(\overline{L}
(T_m^{p_1-(c+j),p_2})_{\phi^{(c+j)}} z^k\overline{z}^k
{L}u,  (T_m^{p_1,p_2})_\phi u)$$ 
\begin{equation}\label{est:ap3}
-\sum_{k'=1}^k \sum_{{1\leq c, 0\leq
j}\atop c+ j\leq p_2}\frac{1}{c!}S^{c+j}_j
(\overline L{z}^{k'}
(T_m^{p_1,p_2-(c+j)})_{\phi^{(c+j)}}
{z}^{k-k'}\overline z^k{L}u, 
(T_m^{p_1,p_2})_\phi u) \end{equation}
$$+\sum_{k'=0}^{k-1}\sum_{{1\leq c, 0\leq j}\atop c+ j\leq
p_1}\frac{1}{c!}S^{c+j}_j(\overline L{z}^k
\overline z^{k'} 
(T_m^{p_1-(c+j),p_1})_{\phi^{(c+j)}}
\overline z^{k-k'}{L}u, 
(T_m^{p_1,p_2})_\phi u)$$ 
$$+\, \sum_{{1\leq c, 0\leq j}\atop c+ j\leq
p_2}\frac{1}{c!}S^{c+j}_j(\overline L{z}^k
\overline z^k
(T_m^{p_1,p_2-(c+j)})_{\phi^{(c+j)}}{L}u,  (T_m^{p_1,p_2})_\phi u) 
$$ 
$$=A_1+A_2+A_3+\sum_{k'=1}^k A_{4,k'}+\sum
_{k'=0}^{k-1}A_{5,k'}+A_6.$$ 

Concerning the critical $L, \overline{L},z^k$ and
 $\overline{z}^k,$ note that in each term above
\begin{itemize}
\item no $L, \overline{L},$ power of $z$ or power
of $\overline{z}$ has been lost,
\item the order among $L, \overline{z}^k, z^k,$ and $
\overline{L}$ is preserved,
\item  letting $|q|=q_1+q_2,$ each term on the right
contains $(T^{q_1,q_2})_{\phi^{(|p|-|q|)}}$ with
$|q|<|p|$ and $|p|-|q|$ derivatives on $\phi,$
\item just as (\ref{est:ap3}) demonstrates the errors
which result in moving
$(T^{p_1,p_2})_\phi$ past the vector fields $L, \olm,
z^k\olm$ and $\oz^kL,$ further such brackets to position
the vector fields so as to make use the {\it a priori}
estimate again will produce similar errors, with
$|q|$ still lower and the `lost' $T$ derivatives
transferred to $\phi,$ 
\item iterating this process, together with a weighted
Schwarz inequality, will produce a sum of terms 
with $q_j\leq \frac{p}{2}$ of the
form 
$$l.c.\,\|\Lambda^{\frac{k-1}{2}}
(T_m^{q_1,q_2})_{\phi^{(|p|-|q|)}}Pu\|_0^2
+s.c.\,\|\Lambda^{-\frac{k-1}{2}}
(T_m^{\frac{p}{2},\frac{p}{2}})_\phi u)\|^2.$$
\item In fact, using larger constants
$\tilde S^{c+j}_j$ subject to the same kind
of bounds,
$|\tilde{S}^{c+j}_j|\leq
\tilde{C}^{c+j},$ we may replace all sums on the right
hand side above by suprema subject to the same range
restrictions on the indices.
\item This use of suprema allows us easily to iterate
everything on the right with easy control on the
constants until either $p_1$ or
$p_2,$ both of which start as $\frac{p}{2},$ drops to
zero, which may happen in two ways - either by
stepwise decrease as on the right hand side above from
successive brackets or by the single term in
Propositions \ref{3.1} and \ref{3.3} which is not
cancelled, the term with all $L$'s or $\ol$'s on one
side or the other in the definition of
$(T_m^{p_1,p_2})_\phi,$ whose principal term is 
$(T_m^{p_1,0})_{\phi^{(p_2)}} (zL)^{p_2}/p_2!$ or its
analogue with $p_1$ reduced to $0.$
\item At this point we no longer have an 
effective localization of powers of $T$ - for
example, brackets with $\ol$ are not corrected. We
proceed anyway, and when we lack a `good' vector field
such as $\ol$ (or of course $\oz^kL$), we create one by
integrating by parts:
\begin{equation}\label{lol}\|Lw\|^2 \leq \|\ol w\|^2 +
|((|z|^{2m})Tw,w)|\end{equation} 
to use up the $L$ and $\ol$ derivatives with the
byproduct of introducing up to half the 
number of new $T$ derivatives. 
\end{itemize}

Overall, then, the strategy has
been, with universal $C$, bounding derivatives of the
Ehrenpreis-type localizing functions 
$\phi^{(r)}$ by $(CN)^r\sim C^rr!$ for $r\leq N\sim p$:
\begin{equation}\label{pto3p/2}
{\|\Lambda^{-\frac{k-1}{2}}T^{p}u\|_{\{\phi \equiv
1\}}\over N^p}\rightarrow
{\|\Lambda^{-\frac{k-1}{2}}(T_m^{\frac{p}{2},\frac{p}{2}})_\phi
u\|\over N^p}\rightarrow\end{equation}
$$\rightarrow 
C^{\frac{p}{4}}{\|\Lambda^{-\frac{k-1}{2}}T^{\frac{3p}{4}}u\|_{supp
\phi}\over N^\frac{3p}{4}}\rightarrow 
C^{\frac{p}{4}}{\|\Lambda^{-\frac{k-1}{2}}T^{\frac{3p}{4}}u\|_{\{
\phi_1\equiv 1\}}\over
N^\frac{3p}{4}}\rightarrow\ldots$$ for suitable $\phi_1
\equiv 1$ on the support of $\phi.$ This will
continue, with a sequence of $\phi_j$ supported
in nested intervals as in \cite{Ta1978},
\cite{Ta1980} until only a negligible fraction
of $p$ is left, namely a bounded number of
derivatives. Since the order in $T^p$ is reduced by
a factor of $3/4$ each time, we will need
$\log_{4/3}p$ such nested open sets. Thus: 
\begin{equation}\label{pto3p/2}{\|\Lambda^{-\frac{k-1}{2}}T^{p}u\|_{\{\phi
\equiv 1\}}\over N^p}\
\lesssim C^p\|\Lambda^{-\frac{k-1}{
2}+3}u\|_{\{\phi_{\log_\frac{4}{3}p}\equiv
1\}}\end{equation}
where the $3$ could be any other small integer. And of
course the whole derivation could have been done at the
$H^s$ level: for any given $s,$ 
\begin{equation}\label{pto3p/2}
{\|\Lambda^{s-\frac{k-1}{2}}T^{p}u\|_{\{\phi \equiv
1\}}\over N^p}\
\lesssim
C^p\|\Lambda^{s-\frac{k-1}{2}+3}u\|
_{\{\phi_{\log_\frac{4}{3}p}\equiv
1\}}\end{equation} which will end the story if this last
norm is known to be finite even as the localizer in
$\tau$ tends to the identity, provided 
that the terms that arise
along the way are all similarly bounded. The most
important of these is of course 
$$|(\Lambda^{s+\frac{k-1}{2}}(T_m^{\frac{p}{2},\frac{p}{2}})_\phi
Pu,
\Lambda^{s-\frac{k-1}{2}}(T_m^{\frac{p}{2},\frac{p}{2}})_\phi
u)|$$ which shows that $P_mu\in H^{s+\frac{k-1}{2}+p}$
implies that $u\in H^{s-\frac{k-1}{2}+p},$ a loss
of 
${k-1}$ derivatives.
 
\begin{remark} The value of $s$ will be chosen so that
we know the norm on the right in (\ref{pto3p/2}) is
finite as the localizer in $\tau$ goes to the identity,
and then $p$ will be chosen so that $P_mu\in 
H^{s+\frac{k-1}{2}+p}$ for that value of $s.$ It follows
that $u\in H^{s-\frac{k-1}{2}+p}.$
\end{remark} 
\begin{remark} For analyticity, one needs to ensure
that as we take $p$ larger and larger, the constants
satisfied by the Ehrenpreis-type localizers are
subject to bounds such that the estimate
(\ref{pto3p/2}) is uniform in $p.$ We have shown this
often before and the arguments are the same here.
\end{remark}

\end{document}